\documentclass[12pt,a4paper,reqno]{amsart}
\usepackage{amsmath}
\usepackage{amsfonts}
\usepackage{amssymb}
\usepackage{amsthm}
\usepackage{enumerate}
\usepackage{centernot}
\usepackage{fullpage}
\usepackage{mathrsfs}
\usepackage{graphicx}
\usepackage{color}

\newcommand{\upperRomannumeral}[1]{\uppercase\expandafter{\romannumeral#1}}

\theoremstyle{plain}
  
  \newtheorem{claim}[]{Claim}
  \newtheorem{proposition}[]{Proposition}
  \newtheorem{lemma}[]{Lemma}
  \newtheorem{theorem}[]{Theorem}

  \newtheorem{remark}[]{Remark}

 \def \J{{\mathbf J}}
 \def \x{{\mathbf x}}
  \def \0{{\mathbf 0}}
 \def \m{{\mathbf m}}
 \def \n{{\mathbf n}}
 \def \k{{\mathbf k}}
 \def\N{{\mathbb N}}

\title[Ursell functions]{Erratum to: Monotonicity of Ursell functions in the Ising model}

\author{Federico Camia}
\address{Division of Science, NYU Abu Dhabi, Saadiyat Island, Abu Dhabi, UAE \& Courant Institute School, New York University, 251 Mercer Street, New York, NY 10012, USA.}
\email{federico.camia@nyu.edu}
\author{Jianping Jiang}
\address{Yau Mathematical Sciences Center, Tsinghua University, Beijing 100084, China.}
\email{jianpingjiang@tsinghua.edu.cn}
\author{Charles M. Newman}
\address{Courant Institute School, New York University, 251 Mercer Street, New York, NY 10012, USA \& NYU-ECNU Institute of Mathematical Sciences at NYU Shanghai, 3663 Zhongshan Road North, Shanghai 200062, China.}
\email{cn3@nyu.edu}

\begin{document}


\maketitle

It was pointed out in the appendix of~\cite{BCS26} that the inequality (43) in ~\cite[Proposition~2]{CJN23}, used in the proof of ~\cite[Theorem~1]{CJN23}, is incorrect in general. We provide a proof of ~\cite[Theorem 1]{CJN23} here; the proof of \cite[Theorem~2]{CJN23} based on ~\cite[Theorem 1]{CJN23} is correct. To make the erratum self-contained, we recall some definitions and Theorem~1 from~\cite{CJN23}.

Let $G=(V,E)$ be a finite graph with vertex set $V$ and edge set $E$. The Ising model on $G$ with \textit{ferromagnetic pair interactions} (or \textit{couplings}) $\J:=(J_{uv})_{uv\in E}$ where $J_{uv}\in[0,\infty)$ for each $uv\in E$ and external field $\mathbf{h}:=(h_v)_{v\in V}$ is defined by the probability measure $\mathbb{P}_{G,\J, \mathbf{h}}$ on $\{-1,+1\}^V$ such that
\begin{equation}\label{eq:Isingdef}
	\mathbb{P}_{G,\J, \mathbf{h}}(\sigma)=\frac{\exp\left[\sum_{uv\in E}J_{uv}\sigma_u\sigma_v+\sum_{u\in V} h_u\sigma_u\right]}{Z_{G,\J, \mathbf{h}}}, \sigma\in\{-1,+1\}^V,
\end{equation}
where $Z_{G,\J, \mathbf{h}}$ is the partition function 
\[Z_{G,\J, \mathbf{h}}:=\sum_{\sigma \in \{-1,+1\}^V}\exp\left[\sum_{uv\in E}J_{uv}\sigma_u\sigma_v+\sum_{u\in V} h_u\sigma_u\right].\]
Let $\langle \cdot \rangle_{G,\J, \mathbf{h}}$ be the expectation with respect to $\mathbb{P}_{G,\J, \mathbf{h}}$. For $k\in \mathbb{N}$ and $j_1,\dots,j_{k} \in V$, the \textit{Ursell function} $u_k(\sigma_{j_1},\dots,\sigma_{j_k})$ for the zero-field Ising measure $\mathbb{P}_{G,\J, \mathbf{0}}$ is defined by
\begin{equation}\label{eq:ursdef}
	u_k(\sigma_{j_1},\dots,\sigma_{j_k}):=\left.\frac{\partial^k}{\partial h_{j_1}\dots\partial h_{j_k}}\ln\left\langle \exp\left[\sum_{u\in V} h_{u}\sigma_{u}\right]\right\rangle_{G,\J, \mathbf{0}}\right|_{\mathbf{h}=\mathbf{0}}.
\end{equation}

The spin-flip symmetry implies that
\begin{equation}
	u_{2k-1}(\sigma_{j_1},\dots,\sigma_{j_{2k-1}})=0, \forall k\in\mathbb{N}, \forall j_1,\dots,j_{2k-1}\in V.
\end{equation}

\begin{theorem}\label{thm:umon}
	Let $G=(V,E)$ be a finite graph and $\J$ be ferromagnetic pair interactions on $G$. Then for any $k\in\mathbb{N}$, any $j_1,\dots,j_{2k}\in V$ ($j_l$'s are not necessarily distinct), and any $u_0v_0\in E$, we have
	\begin{equation}\label{eq:umon}
		(-1)^{k-1}\frac{\partial u_{2k}(\sigma_{j_1},\dots,\sigma_{j_{2k}})}{\partial J_{u_0v_0}}\geq 0.
	\end{equation}
\end{theorem}
\begin{remark}
	When all couplings are zero, $u_{2k}(\sigma_{j_1},\dots,\sigma_{j_{2k}})=0$ if $j_1,\dots,j_{2k}$ contain at least two distinct vertices and $(-1)^{k-1}u_{2k}(\sigma_{u},\dots,\sigma_{u}) \geq 0$ by an explicit computation as in~\eqref{eq:lncos} below; a corollary of Theorem \ref{thm:umon} is that $(-1)^{k-1} u_{2k} \geq 0$~\cite{Shl86}.
\end{remark}

We use the following notation for truncated two-point and three-point functions:
\begin{align*}
	&\langle \sigma_u; \sigma_v \rangle:=\langle \sigma_u \sigma_v \rangle-\langle \sigma_u \rangle \langle \sigma_v \rangle, \\
	& \langle \sigma_u; \sigma_v; \sigma_w \rangle:=\langle \sigma_u \sigma_v \sigma_w \rangle-\langle \sigma_u \sigma_v \rangle \langle \sigma_w \rangle-\langle \sigma_u \sigma_w \rangle \langle \sigma_v \rangle-\langle \sigma_v \sigma_w \rangle \langle \sigma_u \rangle +2\langle \sigma_u \rangle \langle \sigma_v \rangle \langle \sigma_w \rangle.
\end{align*}
For $\x \in \mathbb{R}^V$ with $\|\x\|_1 < \pi/2$, define the characteristic function and cumulant-generating function
\begin{equation}\label{eq:phiFdef}
	\phi_{G,\J}(\x):=\langle \exp[i\sum_{u\in V} x_u\sigma_u]\rangle_{G,\J, \0}\geq \cos\|\x\|_1>0,~F_{G,\J}(\x):=-\ln \phi_{G,\J}(\x).
\end{equation}
Note that
\begin{equation}\label{eq:u1-3}
	\frac{\partial F_{G,\J}(\x)}{\partial x_u}=-i \langle \sigma_u \rangle_{G,\J,i\x}, \frac{\partial^2 F_{G,\J}(\x)}{\partial x_v\partial x_u}= \langle \sigma_u; \sigma_v \rangle_{G,\J,i\x}, \frac{\partial^3 F_{G,\J}(\x)}{\partial x_w\partial x_v \partial x_u}=i \langle \sigma_u; \sigma_v; \sigma_w \rangle_{G,\J,i\x}.
\end{equation}
For  $uv\in E$, define
\begin{equation}\label{eq:Hdef}
	H_{G,\J,uv}(\x):=\frac{\partial F_{G,\J}(\x)}{\partial J_{uv}}=\langle \sigma_u\sigma_v \rangle_{G,\J,\0}-\langle \sigma_u\sigma_v \rangle_{G,\J,i\x}, \  \forall \|\x\|_1<\pi/2.
\end{equation}
Both $F$ and $H$ are real-analytic in $\{ \x \in \mathbb{R}^V: \| \x \|_1<\pi/2\}$. Since
\[(-1)^{k-1}\frac{\partial u_{2k}(\sigma_{j_1},\dots,\sigma_{j_{2k}})}{\partial J_{u_0v_0}}=\left.\frac{\partial}{\partial J_{u_0v_0}}\frac{\partial^{2k}F_{G,\J}(\x)}{\partial _{x_{j_1}}\dots \partial_{x_{j_{2k}}}} \right|_{\x=\0}=\left. \frac{\partial^{2k}H_{G,\J,u_0v_0}(\x)}{\partial _{x_{j_1}}\dots \partial_{x_{j_{2k}}}} \right|_{\x=\0},\]
the proof of Theorem~\ref{thm:umon} boils down to the proof of 
\begin{proposition}\label{prop:Taylor}
	Let $G$ and $\J$ be defined as in Theorem~\ref{thm:umon}. For each $k\in \N$, any $j_1,\dots,j_{2k}\in V$, with arbitrary repetitions, and any $u_0v_0\in E$,
	\[\left. \frac{\partial^{2k}H_{G,\J,u_0v_0}(\x)}{\partial _{x_{j_1}}\dots \partial_{x_{j_{2k}}}} \right|_{\x=\0} \geq 0.\] 
\end{proposition}

We prove Proposition~\ref{prop:Taylor} by induction. We first use the random current representation to build a recursive relation between the second-order partial derivative of $F_{G,\J}$ and $F$ defined on subgraphs of $G$. For $U\subset V$, let $G[U]$ be the induced subgraph with vertex set $U$ and edge set $E(G[U]):=\{uv\in E: u,v \in U\}$. A \textit{current} $\n$ on $G[U]$ is a mapping from $E(G[U])$ to $\N_0:=\N\cup \{0\}$. A \textit{source} of $\n$ is a vertex $u\in U$ at which $\sum_{v\sim u} \n_{uv}$ is odd where $v\sim u$ means $uv \in E(G[U])$; the set of all sources of $\n$ is denoted by $\partial \n$.
\begin{lemma}\label{lem:recursive}
	Let $G$ and $\J$ be defined as in Theorem~\ref{thm:umon}. For any $u_0,v_0\in V$, not necessarily distinct, we have
	\begin{equation}\label{eq:recursive}
		\langle \sigma_{u_0}; \sigma_{v_0} \rangle_{G,\J,i\x}=\sum_{U\supset \{u_0,v_0\}} w_U^{u_0v_0}\exp[2(F_{G,\J}(\x)-F_{G[U^c],\J}(\x))],
	\end{equation}
	where $U^c:=V\setminus U$, $Z_{\emptyset,\J,0}:=1$, $F_{\emptyset,\J}:=0$ and
	\begin{equation}\label{eq:omega_Udef}
		w_U^{u_0v_0}:=\frac{Z^2_{G[U^c],\J,\0}}{Z^2_{G,\J,\0}} 2^{1+|U|}\sum_{\substack{\n \in \N_0^{E(G[U])}, \partial \n=\{u_0\}\Delta \{v_0\}\\ (U,\{e:\n_e>0\}) \emph{ connected}}} \prod_{e\in E(G[U])} \frac{(2J_e)^{\n_e}}{\n_e!} \geq 0,
	\end{equation}
	where $A\Delta B:=(A\setminus B)\cup (B\setminus A)$ is the symmetric difference.
\end{lemma}
\begin{proof}
	Let $\sigma \sim \mathbb{P}_{G,\J,i\x}$ and $\tau \sim \mathbb{P}_{G,\J,i\x}$ be independent. Define
	\[s_u:=\frac{\sigma_u+\tau_u}{2},~t_u:=\frac{\sigma_u-\tau_u}{2},\ \forall u\in V.\]
	Then $(s_u,t_u)\in \Omega:=\{(\pm1,0),(0,\pm1)\}$. The joint Hamiltonian is
	\[\sum_{uv \in E} J_{uv}(\sigma_u\sigma_v+\tau_u\tau_v)+i\sum_{u\in V}(\sigma_u+\tau_u)x_u=2\sum_{uv\in E}J_{uv}(s_us_v+t_ut_v)+2i\sum_{u\in V} x_us_u.\]
	So $(s,t)$ has distribution
	\[P_{G,\J,i\x}(s,t) \propto \exp[2\sum_{uv\in E}J_{uv}(s_us_v+t_ut_v)+2i\sum_{u\in V} x_us_u].\]
	Write $\langle\langle \cdot \rangle \rangle_{G,\J,i\x}$ for the corresponding expectation. Then
	\begin{align}\label{eq:ursell2}
		\langle \sigma_{u_0}; \sigma_{v_0} \rangle_{G,\J,i\x}&=2\langle\langle t_{u_0}t_{v_0} \rangle \rangle_{G,\J,i\x}\nonumber\\
		&=2\frac{\sum_{(s,t)\in \Omega^V}t_{u_0}t_{v_0}\exp[2\sum_{uv\in E}J_{uv}(s_us_v+t_ut_v)]\exp[2i\sum_{u\in V} x_us_u]}{Z^2_{G,\J,i\x}}.
	\end{align}
Applying the Taylor series of $e^z$, we get that the numerator equals
\begin{align}\label{eq:numerator}
	&\sum_{(s,t)\in \Omega^V}t_{u_0}t_{v_0}\prod_{uv \in E} \sum_{\n_{uv}=0}^{\infty}\frac{(2J_{uv}s_us_v)^{\n_{uv}}}{\n_{uv}!}\sum_{\m_{uv}=0}^{\infty} \frac{(2J_{uv}t_ut_v)^{\m_{uv}}}{\m_{uv}!}\exp[2i\sum_{u\in V} x_us_u]\nonumber\\
	=&\sum_{\n,\m \in \N_0^E} w(\n)w(\m)\sum_{(s,t)\in \Omega^V}\prod_{u \in V}s_u^{\sum_{v\sim u}\n_{uv}} t_u^{\sum_{v\sim u}\m_{uv}+\delta_{u_0,u}+\delta_{v_0,u}}\exp[2i x_us_u],
\end{align}
where $\delta_{u,v}$ is the Kronecker delta, and
\[w(\n):=\prod_{uv\in E} \frac{(2J_{uv})^{\n_{uv}}}{\n_{uv}!}.\]
It is easy to see that
\begin{equation}\label{eq:sumcases}
	\sum_{(s_u,t_u)\in \Omega} s_u^{k}t_u^{l}\exp[2ix_us_u]=\begin{cases}
		0, & k\in \N, l\in \N,\\
		0, & k=0, l \text{ odd},\\
		2, & k=0, l \text{ even}>0.
	\end{cases}
\end{equation}
After summing over $(s,t)\in \Omega^V$, every $(\n,\m)$ in \eqref{eq:numerator}  with nonzero contribution satisfies
\begin{equation}\label{eq:mrestriction}
	\sum_{v\sim u} \m_{uv}+\delta_{u_0,u}+\delta_{v_0,u} \text{ is even}, \ \forall u\in V.
\end{equation}
Let $U$ be the vertex set of the connected component of $(V,\{e\in E: \m_e>0\})$ which contains $u_0$ and $v_0$; if $u_0\neq v_0$, then \eqref{eq:mrestriction} implies that there is a path from $u_0$ to $v_0$ in the graph with vertex set $U$ and edge set $\{e=uv: u, v\in U, \m_{e}>0\}$.

By the first case in \eqref{eq:sumcases}, all non-zero terms in \eqref{eq:numerator} satisfy $\{u,v \in V: e=uv\in E, \n_{e}>0\} \cap U=\emptyset$. The definition of $U$ implies that $\m_{uv}=0$ for all $u\in U$ and $v\in U^c$. Therefore, non-zero terms in \eqref{eq:numerator} satisfy
\[\n_{uv}+\m_{uv}=0, \ \forall uv\in E \text{ with }u\in U, v\in U^c.\]
These and the third case in \eqref{eq:sumcases} imply that~\eqref{eq:numerator} equals
\begin{align*}
	&\sum_{U\supset\{u_0,v_0\}}\sum_{\n,\m \in \N_0^{E(G[U^c])}}w(\n)w(\m)\sum_{(s,t)\in \Omega^{U^c}}\prod_{u\in U^c}\left(s_u^{\sum_{v\sim u} \n_{uv}}t_u^{\sum_{v\sim u}\m_{uv}}\exp[2i x_us_u]\right)\\
	&\qquad\times2^{|U|}\sum_{\substack{\k \in \N_0^{E(G[U])},\partial \k=\{u_0\}\Delta \{v_0\}\\ (U,\{e:\k_e>0\}) \text{ connected}}}w(\k)\\
	=&\sum_{U\supset\{u_0,v_0\}} Z^2_{G[U^c],\J,i\x}2^{|U|}\sum_{\substack{\k\in \N_0^{E(G[U])},\partial \k=\{u_0\}\Delta \{v_0\}\\ (U,\{e:\k_e>0\}) \text{ connected}}}w(\k).
\end{align*}
Plugging this into \eqref{eq:ursell2}, and then using \eqref{eq:omega_Udef}, we get
\begin{align*}
	\langle \sigma_{u_0}; \sigma_{v_0} \rangle_{G,\J,i\x}=\sum_{U\supset\{u_0,v_0\}} w_U^{u_0v_0}\frac{Z^2_{G[U^c],\J,i\x}}{Z^2_{G,\J,i\x}}\frac{Z^2_{G,\J,\0}} {Z^2_{G[U^c],\J,\0}}=\sum_{U\supset\{u_0,v_0\}} w_U^{u_0v_0}\frac{\phi^2_{G[U^c],\J}(\x)}{\phi^2_{G,\J}(\x)}.
\end{align*}
This and the definition of $F$ in \eqref{eq:phiFdef} complete the proof of the lemma.
\end{proof}
\begin{proof}[Proof of Proposition~\ref{prop:Taylor}]
	For a real-analytic function $f$ defined near $\0$ and $n\in \N_0$, we say $f$ has \textit{property $T_n$} if the Taylor expansion of $f$ at $\0$ has multivariate Taylor coefficients nonnegative through total degree $n$:
	\[\left.\frac{\partial^kf(\x)}{\partial x^{k_1}_1\dots\partial x^{k_{|V|}}_{|V|}} \right |_{\x=\0} \geq 0, \ \forall 0\leq k \leq n \text{ with } k_1+\dots+k_{|V|}=k, k_j\in \N_0.\]
	We use induction on even $n$ to prove: for all finite graphs $G=(V,E)$, all nonnegative finite couplings $\J$, all $uv \in E$, $H_{G,\J,uv}$ defined in~\eqref{eq:Hdef} has property $T_n$ for each even $n\in \N_0$. Note that the Taylor coefficients of $H_{G,\J,uv}$ with odd total degree are zero since $H_{G,\J,uv}(\x)=H_{G,\J,uv}(-\x)$. The base case is trivial since $H_{G,\J,uv}(\0)=0$. We now assume that for all such $G,\J$ and $uv$, $H_{G,\J,uv}$ has property $T_{n-2}$ for some even $n\geq 2$. Our goal is to show that $H_{G,\J,uv}$ has property $T_{n}$; and this is achieved via the following three claims.
	\begin{claim}\label{claim1}
		For all finite graphs $G=(V,E)$, all nonnegative finite $\J$, all $U\subset V$, $F_{G,\J}-F_{G[U^c],\J}$ has property $T_{n-2}$.
	\end{claim}
	For $\lambda \in [0,\infty)$, define $\J^{\lambda}$ by
	\begin{equation*}
		J^{\lambda}_e:=\begin{cases}
			J_e,& e\in E(G[U^c]),\\
			\lambda J_e,& e \cap U \neq \emptyset,
		\end{cases}
	\end{equation*}
	where we abuse the notation $e \cap U \neq \emptyset$ for $e\notin E(G[U^c])$.
	Then the definition of $H$ gives
	\[\frac{\partial \left(F_{G,\J^{\lambda}}-F_{G[U^c],\J^{\lambda}}\right)}{\partial \lambda}=\sum_{e: e\cap U\neq \emptyset} J_e H_{G,\J^{\lambda},e}.\]
	Since 
	\[F_{G,\J^0}(\x)-F_{G[U^c],\J^0}(\x)=-\ln \prod_{u\in U}\phi_{\{u\},\0}(x_u)=-\sum_{u\in U} \ln \cos(x_u),\]
	we have
	\[F_{G,\J}(\x)-F_{G[U^c],\J}(\x)=-\sum_{u\in U} \ln \cos(x_u)+\sum_{e: e\cap U\neq \emptyset} J_e \int_0^1 H_{G,\J^{\lambda},e}(\x) d \lambda.\]
	All the Taylor coefficients of $-\ln \cos(x_u)$ are nonnegative because for all $|x|<\pi/2$,
	\begin{equation}\label{eq:lncos}
		\cos x=\prod_{k=0}^{\infty}\left(1-\frac{x^2}{a_k^2}\right) \text{ with }a_k:=\frac{(2k+1)\pi}{2},\  -\ln \cos(x)=\sum_{n=1}^{\infty} \left(\sum_{k=0}^{\infty} \frac{1}{na_k^{2n}}\right) x^{2n}.
	\end{equation}
	So the induction hypothesis proves Claim~\ref{claim1}.
	\begin{claim}\label{claim2}
		For all finite graphs $G$, all nonnegative finite $\J$, $F_{G,\J}$ has property $T_n$.
	\end{claim}
	This follows from~\eqref{eq:u1-3},  Lemma~\ref{lem:recursive} and Claim~\ref{claim1}. Indeed, up to a positive multiplicative constant, Taylor coefficients of $F_{G,\J}$ with total degree in $[2, n]$ can be written as Taylor coefficients of $\partial^2F_{G,\J}(\x)/\partial x_u\partial x_v$ for some $u,v\in V$ with total degree $\leq n-2$; the constant coefficient is $F_{G,\J}(\0)=0$ and linear coefficients vanish by spin-flip symmetry.
	
	For each $w\in V$, differentiate \eqref{eq:recursive} with respect to $x_w$ and apply \eqref{eq:u1-3} to obtain
	\begin{align}\label{eq:recursive1}
		&-i\left[2\langle \sigma_{u_0}; \sigma_{v_0} \rangle_{G,\J,i\x}\langle \sigma_w \rangle_{G,\J,i\x}+\langle \sigma_{u_0};\sigma_{v_0};\sigma_w\rangle_{G,\J,i\x}\right]\nonumber\\
		=&2\sum_{U\supset \{u_0,v_0\}} w_U^{u_0v_0}\frac{\partial F_{G[U^c],\J}(\x)}{\partial x_w}\exp[2(F_{G,\J}(\x)-F_{G[U^c],\J}(\x))].
	\end{align}
	\begin{claim}\label{claim3}
		For all finite $G=(V,E)$, all nonnegative finite $\J$, all $u_0,v_0,w\in V$, 
		\[-i\left[2\langle \sigma_{u_0}; \sigma_{v_0} \rangle_{G,\J,i\x}\langle \sigma_w \rangle_{G,\J,i\x}+\langle \sigma_{u_0};\sigma_{v_0};\sigma_w\rangle_{G,\J,i\x}\right] \emph{ has property }T_{n-1}.\]
	\end{claim}
	Indeed, in~\eqref{eq:recursive1}, a Taylor coefficient of total degree $\leq n-1$ in the product
	\[\frac{\partial F_{G[U^c],\J}}{\partial x_w}\exp[2(F_{G,\J}(\x)-F_{G[U^c],\J}(\x))]\]
	is nonzero only if one uses Taylor coefficients of the exponential factor of total degree $\leq n-2$ since $\partial F_{G[U^c],\J}(\x)/\partial x_w|_{\x=\0}=0$ by the first equality in \eqref{eq:u1-3} and spin-flip symmetry; Claim~\ref{claim3} then follows from \eqref{eq:recursive1} and Claims~\ref{claim1} and \ref{claim2}.
	
	For each $w\in V$, by \eqref{eq:Hdef}, we have
	\begin{align*}
		\frac{\partial H_{G,\J,uv}(\x)}{\partial x_w}&=-i\langle \sigma_u\sigma_v;\sigma_w\rangle_{G,\J,i\x}=\frac{-i}{2}\left[2\langle \sigma_{u}; \sigma_{w} \rangle_{G,\J,i\x}\langle \sigma_v \rangle_{G,\J,i\x}+\langle \sigma_{u};\sigma_{v};\sigma_w\rangle_{G,\J,i\x}\right]\\
		&\quad +\frac{-i}{2}\left[2\langle \sigma_{v}; \sigma_{w} \rangle_{G,\J,i\x}\langle \sigma_u \rangle_{G,\J,i\x}+\langle \sigma_{u};\sigma_{v};\sigma_w\rangle_{G,\J,i\x}\right].
	\end{align*}
Applying Claim~\ref{claim3} to each bracket together with its factor $-i$, we see that $\partial H_{G,\J,uv}(\x)/\partial x_w$ has property $T_{n-1}$. Therefore, $H_{G,\J,uv}$ has property $T_n$. This completes the induction and proves Proposition~\ref{prop:Taylor}.
\end{proof}


\subsection*{AI statement}
The proof strategy was developed by the authors, and the entire manuscript was written by the authors. AI tools were used for review and proofreading.
\bibliographystyle{abbrv}
\bibliography{reference}

\end{document}